\input amstex
\documentstyle{amsppt}
\document
\magnification=1100
\NoBlackBoxes

\centerline{I.M.GELFAND AND HIS SEMINAR - A PRESENCE   }

\bigskip
\centerline{A.~Beilinson}

\centerline{University of Chicago}
 
\bigskip

\hfill{\it ...How nice to be like a fool} 

\hfill{\it for then one's Way is grand beyond measure}

\medskip

\hfill{From a poem of Tainin Kokusen}

\hfill{given to his student Ry\~okan Taigu, 1790}

\bigskip

\medskip

The mathematical seminar of Israel Moiseevich Gelfand  started each year in the beginning of September and ended  in the spring when IM would observe  that   ``rivulets (of melting snow) are beginning to flow." The sessions   were on Mondays in the  big auditorium on the 14th floor of  Moscow University's main building, and each  consisted of two parts: a  preseminar
that began at 6pm, and  the seminar proper,  which began with IM's arrival
at around 7pm  and ended at  10pm when a cleaning lady entered the room to announce her departure (at which time  the floor was to be locked,  and those  wishing to spend the night at home had to hurry down). During the  preseminar  dozens of people  congregated near the auditorium entrance, chatting  and exchanging books and texts of all kinds.\footnote{Senya Gindikin: ``IM considered these  preseminar discussions   to be very important. However, he was   pathologically disorganized and could not get anywhere on time even if he wished to (e.g., for a meeting with important people)."  } 
  The seminar typically began with IM telling some anecdotes and mathematical news,  after which would come a talk by an invited speaker.\footnote{For Misha Shubin's notes of the talks, see http://www.mccme.ru/gelfand/notes/.} Often  there was not enough time to finish, and the talk continued  serially,  each time beginning from   scratch and covering about half of the material from the week before, the speaker gradually fading away and being replaced by a student assigned by IM to explain what the talk was, or should have been, about. Any speaker  deemed not to have understood the subject,  or to have explained it badly (or if the writing was too small and the voice not clear) was harshly reprimanded.\footnote{At times the scene resembled the   koan about Nansen and the cat,  see http://en.wikipedia. org/wiki/Nanquan\_Puyuan, with no J\~osh\~u in sight. }

\medskip

The seminar started in 1943; I saw its later years,  which coincided with the late period of the Soviet Union. After Stalin's death the edifice of the state shrank into itself, and  free space teemed with  life.   The ideology had lost its fulcrum, 
 the show of democracy was simple  (a single candidate to vote for, not two equally unpalatable ones),
newspapers  were mostly used as toilet tissue. The remaining  taboos  were private commerce\footnote{One exception was Bird's Market in Moscow,  where on the weekends all kinds of animals were sold. Once, when  I visited it with Don Zagier,  a bearded fellow in sheepskins tried to sell him a   snow-white goose. The fellow said that he could  see  Don to be a true gentleman -  otherwise he would not offer him the beauty. He said that she would be Don's best friend,  going with him everywhere, and sharing his bath. The discussion was in French.  } and entrepreneurship,
and political activity outside the Party's womb.  Many people shared 
the attitude of Pushkin's poem ``From Pindemonte"\footnote{For Nabokov's translation, see  
https://ireaddeadpeople.wordpress.com/2014/11/06/alexander-pushkin-to-stroll-in-ones-own-wake/. }  and viewed all matters political as not interesting anyway. The market in a modern sense, this incessant gavage of unneeded things, did not exist. One could quit the tarmac road to look for 
one's own trail into the woods. If the trail happened to be mathematics, it would surely meet with  IM's seminar. 

There was a distinct inner music.\footnote{Perhaps not unlike  that of another closed country - Japan of the  late Edo era. The mores were also not altogether different: e.g., the chief nuclear scientist who dealt with  Chernobyl's   aftermath  killed himself, probably,  as an apology for his involvement in the  nuclear industry (his superiors  practiced the Fukushima era ethics).} 
The air was thin and transparent. 
One could hear the sound of one's breathing, of snowflakes falling,  of hoarfrost's brush decorating the windowpane. Old villages still existed  within Moscow limits, such as wonderful  Dyakovo,   its empty church over an ancient cemetery on a high  scarp above   Moscow River,   wooden houses  edged by deep ravines, and vast   apple gardens where nightingales sang.\footnote{Dyakovo was eliminated in the 80s:   first   the graves in the cemetery were dug out,  then, in a while, the houses were demolished and burned, a single one surviving   for over a year.} Poetry was by far  more real than  social ranks  -  
  poems were rewritten by hand and learned by heart.\footnote{ Two of my friends   knew by heart all of  Mandelstam's poems.    Cf.~"An evening of Russian poetry" by Nabokov,  http://www.sapov.ru/novoe/n00-39.htm. }

\medskip

I was brought to the seminar and introduced to IM by Alesha Parshin in the fall of 1972; I was then a senior of the  2nd mathematical school (IM taught there years earlier).  The precious feeling of being a fool and despite that, or rather thanks to that,  being  in balance with life's flow, akin to running on the cracking ice of a river, goes back to these  days.

After failing the entrance exams to the math department of Moscow University\footnote{The lords of the Moscow mathematical establishment  kept it clean from  anything  Jewish.}
 I found myself in the lovely Pedagogical Institute.  This was a benison - going to classes or playing truant in the morning, going to  mathematical seminars  or taking a train to walk in the woods\footnote{For a year I was seeing the woods almost every day. }  later in the day; and there were wonderful friends. After a while I managed to be transferred to the University. The mood there was more sombre, but with no desire for higher grades,
one could  skip  all the  ideological classes\footnote{Officially one had to know the whole contents of the course in order to pass. But the teachers, with the help of the Komsomol leader of the group, revealed on the eve of the exams  to each student  the exact question he/she would be asked.}  
 to retain the good measure of  idleness and freedom that are  so necessary for doing math.

Accidentally, my first result to be published was close to the one  found at the   same time (the end of 1977) by IM  with Osya Bernstein and Serezha Gelfand. IM gave a talk about his work,
 mentioning  that I had obtained a similar theorem. After the talk I approached IM, and he at once   ordered me to leave Yuri Ivanovich Manin,  who was my supervisor,  and be  his student. The col\'ee was   
 violent. I refused.   When I told YuI about the  
 accolade, he said this had happened to many, e.g., to himself and to Shafarevich. Thereafter  I stayed in an outer orbit of IM's influence, and our relationship  was excellent. 

 After the graduation I got a job in a mathematical laboratory at Moscow Cardiological Center; to that end,  the noble Vladimir Mikhailovich Alexeev, who was the head of the laboratory, came  to the job committee soon after undergoing major cancer surgery. VM died in December of 1980. The new head of the lab, disagreeable on the matter of skiving, was keen on getting rid of me. After IM learned about the situation, he talked to the head of biological sector of the Center;  I was transferred there and left to my own devices. The sinecure   was   better than a   graduate school.

\medskip

In the early 1970s 
the high winds of the Cold War\footnote{Its sole cause was (is) the incompatibility of plutocracy and autocracy; the rest of the US/SU discordances were  red herrings (or, if the reader prefers, forget-me-nots of a Kozma Proutkoff fable that was often cited by IM; see  http://www.math.uchicago.edu/$\sim$mitya/langlands/nezabudki.html  for an English translation).    }   brought permission for Soviet Jews to emigrate, and many signed up for what, with hindsight, was
a verification of the universality of Griboyedov's quip that the place where it is better  for us to be is where we are not.\footnote{These departures, the simulations of dreams,  had little in common with the high quest of crossing the  SU border (in either direction) by one's own free will and with no external purpose,  as in Nabokov's ``Glory" or as  done  by Slava Kurilov, see his   book   ``Alone in the ocean",   http://rozamira.org/lib/names/k/kurilov\_s/kurilov.html (in Russian).} The separation from friends was deemed to be permanent (the imminent demise of the SU was anticipated then no more than that of the US is now).     Dima Kazhdan,  Ilya Iosifovich Piatetski-Shapiro, and Osya Bernstein, with whom we were happily doing math  for his last half year in Moscow, were among those who left. No one at the seminar could replace them.

\medskip

 IM loved  playing with people (with him mischief was never far away).\footnote{Spencer Bloch:  ``I am sure I told you my Gelfand story when he came to Paris and was to meet
with Serre. He was staying at Ormaille and the people at IHES needed someone
to escort him to Paris. I was elected. I suggested we take a train with
plenty of time to spare so we would not inconvenience the great Serre. Of
course, I did not fully grasp the subtle thinking process of my charge.
Suffice it to say that not inconveniencing Serre was rather low on the totem
pole of Gelfand's priorities. I arrived at his apartment and he announced
that he would instruct me on the Russian technique for making tea. So, of
course, we missed the train. But I said no matter, there would be another
train along in 20 minutes. But no, Gelfand said that errors had occurred
during the making of the tea, and nothing would do except to return to his
apartment and make more tea; which we did. So, of course, we missed the next
train. And, as was clearly the intent from the beginning, the great Serre
was made to wait for the great Gelfand."} A common  way to engage someone   was to explore his feeling of self-importance.  
 IM rarely lost the game; if this happened (which meant that the opponent was more unpredictable than IM himself), he was furious, but the winner got his  respect and, perchance, even love.  For example, IM could ask you to wait and then disappear for a very long time.\footnote{ Senya Gindikin: ``I would think that it was more complicated. IM felt no obligations and at every moment did only what he wished to do at that moment. I don't think he did anything intentionally, he could be distracted for a long while. I have a big personal experience here." }\footnote{They say that once, on his way  to a meeting with  President of the Academy of Science, IM stopped to exchange pleasantries with a nice cleaning lady; he did not reach then the destination. }   A cheap win was to leave after an hour. A master stroke would be different.   According to legend, 
 when IM returned to his office after several hours to see how Misha  Tsetlin was doing, he found Misha fast asleep on the sofa.\footnote{Misha Tsetlin, who for IM was what, probably,  J\~osh\~u was for Nansen, died in 1966. About their research in physiology, see sect.~3.1 of M.~Latash's book ``Synergy", Oxford University Press, 2008, http://books.google.ru/books?id=Z45Oj8yCQMIC\&pg=PA53. 
 See also  V.V.Ivanov's article about Tsetlin, http://historyofcomputing.tripod.com/essays/CETLINM.HTM (in Russian).}

IM appreciated life.\footnote{ And, just maybe, he admired its beauty to the point where even ugly human deeds do not blot the clarity of vision. I believe that
 wild animals are  not afraid of   humans who are able to   
 participate that much in the joy of being.}  Although IM was a very social person,  he paid no outward respect to  problems caused by a lack of  inner happiness (as a result he was often perceived as  rude).\footnote{On the other hand, 
IM did care a lot when the problem was real: e.g., his help was crucial for saving the  lives of Sasha Zamolodchikov and the son of Tolya Kushnirenko after terrible accidents.} He did 
what  interested him    for its own sake,   and  not as a part of any grand project.\footnote{IM often said that he abandoned research whenever its subject became too popular.}  Running  seminars (the mathematical and biological\footnote{Volodya Gelfand: ``IM knew no  biology, but was always able to identify true experts to talk to, and these discussions were often very beneficial for the biologists as well."}\footnote{IM was fascinated by biology for  the mystery that you do not know even how to think about is so immediate there.} ones, and, starting from 1986,  the one on informatics)  was always interesting.

Then there was   the work with physicians, a long attempt to find out how a doctor diagnoses  heart disease. While the attempt itself ended in failure,\footnote{ Perhaps at the start IM did not  recognize medicine as an art (for him, a non-mathematician's project to uncover the way a mathematician proves theorems would be laughable). The work on a simpler problem of  diagnosis of meningitis was  successful.} it included several 
top notch physicians  who brought a distinctly new dimension to the life of those around IM.  
I came then to know three doctors,   true masters, who found it impossible to accept any payment for their help.\footnote{ The payment for a cab to bring them, after the  workday at hospital, to the patient's home included.} I learned  that such an attitude is utterly natural and, in truth, a doctor cannot behave any differently.\footnote{A simple criterion to check if a given human society is not dead at its core is the presence in it  of such physicians.}

\medskip

  IM emphasized the  importance of  decency.\footnote{Dima Leshchiner: ``I recall  his favorite saying: ``People do not have shortcomings, but only peculiarities."  It seems to me this has to do with what  ``decency" meant in his understanding, namely,   that ``decency" is the quality of an action, not of a person."  }   Its two realizations central, in my opinion, to IM's life were severing,   after the work on the bomb, his ties with the military (late 1950s)\footnote{IM once told me that, back then, he was offered to be the head of any institution of his choice (say, Institute of Applied Mathematics that dealt with military projects), and he refused. }\footnote{Senya Gindikin: ``I am not sure if anyone knows how and why he stopped the military activity. To what degree this was initiated by himself. He was extremely cautious. He received a closed Lenin prize around 1960."} 
and his becoming  a vegan (mid 1990s).\footnote{See IM's interview for VITA,  http://israelmgelfand.com/talks/vita.html.  In   earlier years IM coauthored a series of works on neurophysiology  based on grisly  experiments on cats.}  Both have to do with 
 overcoming the habit of what is usually called thinking objectively, i.e., paying no regard to violence directed towards others.\footnote{The trite  lament that the cause of the sadness of today's world is that development of  technology has overrun our moral development misses the point - for there is no moral development. Common decency now is the same as it was thousands years ago, and it works well if  applied (and if  those who apply it are not killed). For example, having taken it as a religious principle  (see https://en.wikipedia.org/wiki/Jainism\#Doctrine), the Jains built a reasonable, i.e., nondestructive, society (maybe the only one still in existence). Their  cousins in the West, the  Good People  (referred to as Cathars, ``catlovers", by the adversaries), were eliminated in a feat of what now is called ``globalization".} Arguably, without the first decision the world around IM would have been much less colorful and the seminar  quite different.  Becoming vegetarian is probably  no less  essential.  It may  loosen  knots tied   hard in one's mind,   bringing  back an ability  to see  many things as obvious and simple.

One difference between IM's seminar and other great mathematical seminars was its openness: the talks were not aimed at explaining any distinctive   subject, nor were they connected to IM's current work, but   rather  these were stories that might contain a call from the future. This was  in tune with the next feeling:
We are used to seeing science's accomplishments as being fundamental. Over time   the magical picture switches, and we realize that, in fact, we know almost nothing about the world, and science  merely  attempts  to  hide the vast openness. But we are able to wonder  and   take in new things, and to feel gratitude,
only due to the wind that blows through us.

IM often said that he does not consider himself to be clever.\footnote{``You should not explain to me that I am an idiot: I know this myself."  Instructing Oleg Ogievetsky's mother on how to talk to physicians, IM asserted  ``No one can revoke your inherent right to be a fool."} 
A fool's way to see things differs from that of a clever person like peripheral vision differs from central vision. At every moment there are infinitely many possible directions to look at and to choose. A fool retains awareness of that; a clever person   moves successfully  in one or two directions  while forgetting  completely  the remaining infinity of dimensions. A   new understanding or a fresh poem  starts with a  tiny movement  into an unknown dimension, which is the inimitable act of a fool. 

\medskip

Modern mathematics is a unique thrust of conceptual thought: once the right concept 
  (a mathematical structure) and a language to deal with it are found,   a whole new world unfolds.\footnote{ A related fact is that in mathematics, unlike elsewhere, wrong notions
die off easily.  Our capacity for understanding  is hampered, foremost, by  the inability  to dispel  false concepts.}  So for a mathematician it is very tempting to  search for an adequate language as a key for understanding non-mathematical subjects, e.g., biology. This vision was dear to IM.\footnote{See his Kyoto lecture, http://israelmgelfand.com/talks/kyoto.html, and a birthday party talk, http://www.math.harvard.edu/conferences/unityofmath\_2003/talks/gelfand-royal-talk.html. } One reason why it has not been realized 
might be the following:

Science invariably considers  reality as if from outside, the objects of study clearly distinct from the observer. But mathematical structures are part of the true reality  
that  can be seen only from inside, the object of study being inseparable from the activity of our brain. It might be that adequate languages are peculiar to exactly this type of seeing. 
For example,  except on the most superficial level, science   is blank about the ways animals  interact with the world.  The animal's vision   should be so  wonderfully  different to the human's  that being privy to it  might drastically change our understanding of what reality is. It is in such a quest that an  adequate language could be ignited. Which is a mere foolish dream as long as we persist
 in positioning ourselves as  separated from other living beings and above them - to the degree of imagining  that the earth, the animals, and the trees  can be our property.
 Incidentally, this same delusion   underlies the drive for the destruction of the planet (which has accelerated  so much  since I saw IM the last time).

\medskip

As I am writing these lines, it is spring, the season when the past does not seem to be that 
impossibly separated from the future. Great seminars have something of faerie horses in their nature.
 Bayard is said to be still living somewhere since his escape into the heart of the wild forest of Ardennes.

\medskip

\centerline{********}

This essay would not exist without many walks and discussions with Jesse Ball, Spencer Bloch, Irene(!) and Nicodemus Beilinson, Volodya Drinfeld, Dennis Gaitsgory, Anyuta and Volodya Gelfand, Senya Gindikin, Dima Kazhdan, Dima Leshchiner, Yuri Manin, Oleg Ogievetsky, and Eric Shutt, the request of Slava Gerovitch to write it down,
and  the interest  and help of  Allyn Jackson. My deep gratitude to them.

\enddocument
\end